\documentclass[AMA,Times1COL]{WileyNJD-v5} 
\usepackage{subcaption}
\usepackage{multirow}
\usepackage{caption}
\usepackage{booktabs}
\usepackage{xcolor}
\usepackage{bm}

\lstdefinelanguage{Julia}%
{morekeywords={abstract,break,case,catch,const,continue,do,else,elseif,%
		end,export,false,for,function,immutable,import,importall,if,in,%
		macro,module,otherwise,quote,return,switch,true,try,type,typealias,%
		using,while},%
	sensitive=true,%
	alsoother={$},%
	morecomment=[l]\#,%
	morecomment=[n]{\#=}{=\#},%
	morestring=[s]{"}{"},%
	morestring=[m]{'}{'},%
}[keywords,comments,strings]%

\newcommand\emachine{\epsilon_{\text{machine}}}
\newcommand\bH{\bm H}
\newcommand\bI{\bm I}
\newcommand\bF{\bm F}
\newcommand\bU{\bm U}
\newcommand\bN{\bm N}
\newcommand\bv{\bm v}
\newcommand\bn{\bm n}
\newcommand\bR{\bm R}
\newcommand\bS{\bm S}
\newcommand\btau{\bm \tau}
\newcommand\bC{\bm C}
\newcommand\bb{\bm b}
\newcommand\bE{\bm E}
\newcommand\be{\bm e}
\newcommand\bulk{k}
\newcommand\firstlame{\mathsf{\lambda}}
\newcommand\tcolon{\!:\!}
\newcommand\Det[1]{\lvert #1 \rvert}
\DeclareMathOperator\trace{trace}
\DeclareMathOperator\fl{fl}
\newcommand\expm{\operatorname{\mathtt{expm1}}}
\newcommand\logp{\operatorname{\mathtt{log1p}}}
\newcommand\logpmx{\operatorname{\mathtt{log1p\_minus\_x}}}
\newcommand\Jm{\mathtt{J_{-1}}}
\usepackage{orcidlink}

\newcommand{\cmark}{\checkmark}%
\newcommand{\xmark}{--}%

\articletype{Research Article}%

\received{Date Month Year}
\revised{Date Month Year}
\accepted{Date Month Year}
\journal{International Journal for Numerical Methods in Engineering}
\volume{00}
\copyyear{2024}
\startpage{1}

\begin{document}

\title{Stable numerics for finite-strain elasticity}

\author[1]{Rezgar Shakeri} 
\author[2]{Leila Ghaffari} 
\author[2]{Jeremy L. Thompson} 
\author[2]{Jed Brown*} 

\authormark{SHAKERI \textsc{et al}}
\titlemark{STABLE NUMERICS FOR FINITE-STRAIN ELASTICITY}

\address[1]{\orgdiv{Department of Civil, Environmental, and Architectural Engineering}, \orgname{University of Colorado Boulder}, \orgaddress{\state{CO}, \country{United States}}}

\address[2]{\orgdiv{Department of Computer Science}, \orgname{University of Colorado Boulder}, \orgaddress{\state{CO}, \country{United States}}}

\corres{*Jed Brown \orcidlink{0000-0002-9945-0639}, \email{jed.brown@colorado.edu}}

\fundingInfo{U.S. Department of Energy, Office of Science, Office of Advanced Scientific Computing Research, applied mathematics program and by the National Nuclear Security Administration, Predictive Science Academic Alliance Program (PSAAP) under Award Number DE-NA0003962.}



\abstract[Abstract]{A backward stable numerical calculation of a function with condition number $\kappa$ will have a relative accuracy of $\kappa\emachine$.
	Standard formulations and software implementations of finite-strain elastic materials models make use of the deformation gradient $\bm F = I + \partial \bm u/\partial \bm X$ and Cauchy-Green tensors.
	These formulations are not numerically stable, leading to loss of several digits of accuracy when used in the small strain regime, and often precluding the use of single precision floating point arithmetic.
	We trace the source of this instability to specific points of numerical cancellation, interpretable as ill-conditioned steps.
	We show how to compute various strain measures in a stable way and how to transform common constitutive models to their stable representations, formulated in either initial or current configuration.
	The stable formulations all provide accuracy of order $\emachine$.
	In many cases, the stable formulations have elegant representations in terms of appropriate strain measures and offer geometric intuition that is lacking in their standard representation.
	We show that algorithmic differentiation can stably compute stresses so long as the strain energy is expressed stably, and give principles for stable computation that can be applied to inelastic materials.}



\maketitle



\section{Introduction}\label{introduction}

Errors in computational mechanics are attributable to three sources: continuum model specification (materials, geometry, boundary conditions), discretization (finite elements), and numerical.
When working in double precision with direct solvers, the first two typically dominate numerical errors and stable numerics are overlooked beyond linear algebra.
Meanwhile, single precision is widely considered to be insufficient for finite-strain implicit analysis and practitioners opt for distinct small-strain formulations due to a combination of instability and perceived cost of finite-strain formulations in small-strain regimes.
This shifts a cognitive burden to the practitioner who must confirm that the small-strain formulations are valid, and is problematic for high-contrast materials in which finite strains and infinitesimal strains are present within the same analysis.
In this paper, we demonstrate the instability in standard formulations for hyperelasticity and present intuitive (and mathematically equivalent) reformulations that are stable, enabling finite-strain analysis at all strains and opening the door for reduced precision analysis.

In floating point arithmetic, $\emachine = \sup_{x} \lvert \fl(x) - x \rvert / \lvert x \rvert$ is the maximum error incurred rounding a real number $x$ to its nearest floating point representation $\fl(x)$. (We assume $x$ is within the exponent range and will not discuss overflow/underflow/denormals.)
Typical values of $\emachine$ are $2^{-53} \approx 10^{-16}$ for IEEE-754 double precision and $2^{-24} \approx 6 \cdot 10^{-8}$ for single precision.
Elementary math operators $\circledast$ (standing for addition, subtraction, multiplication, or division) and special functions behave as ``exact arithmetic, correctly rounded'', $x \circledast y = \fl(x \ast y)$ \cite{trefethen1997nla}, thus guaranteeing $\lvert x \circledast y - x \ast y\rvert / \lvert x \ast y \rvert \le \emachine$.

\begin{figure} [h]
	\includegraphics[width=\textwidth]{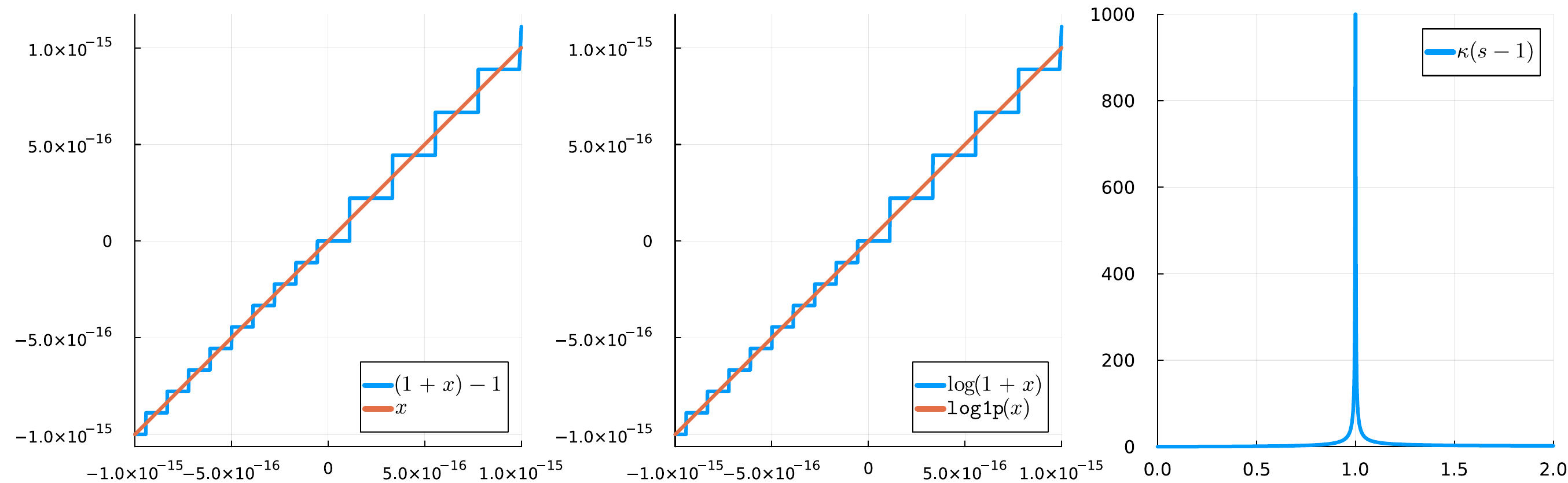}
	\caption{The stair-step patterns illustrate numerical instability evaluating in $(1 + x)-1$ (left) and $\log(1+x)$ (center). When evaluating $(1 + x) - 1$, a small relative error is incurred in the well-conditioned operation $s := 1 + x$, and that error becomes a large relative error because $s - 1$ has unbounded condition number as $s\to 1$ (right), despite this operation being computed exactly in floating point arithmetic.}
	\label{fig:log1p(x)-1-plus-x-m1}
\end{figure}

With such strong guarantees from elementary operations, we might hope that $(x \circledast y) \circledast z$ is also accurate to $\emachine$, but also, this is not true.
An illustrative example is computing $(x + 1) - 1$. For $x = 10^{-16}$ in double precision, $x \oplus 1 = 1$, thus $(x \oplus 1) \ominus 1 = 0$ has a relative error of $1$, which is much larger than $\emachine$. For larger values of $x$, we observe the stair-step effect in \autoref{fig:log1p(x)-1-plus-x-m1}.
The first operation incurred a relative error smaller than $\emachine$ and the second operation was exact, so where can we place blame for the catastrophic cancellation error?
To shed light on this, we consider the condition number of a differentiable function $f(x)$, defined as
\begin{equation}\label{eq:def-cond-number}
	\kappa_{f}(x) = \left\lvert \frac{\partial f}{\partial x} \right\rvert \frac{\lvert x \rvert}{\lvert f(x) \rvert} .
\end{equation}
The first operation $x + 1$ has a condition number of 1 while the second has enormous condition number.
An numerical algorithm $\tilde f$ for evaluating the continuum function $f$ is called \emph{backward stable} if it computes the exact answer to a problem with almost the given inputs, i.e., $\tilde f(x) = f(x')$ for some nearby $x'$ satisfying $\lvert x' - x \rvert/\lvert x \rvert < c \emachine$ for a small constant $c$.
Backward stable algorithms satisfy the forward error bound \cite{trefethen1997nla}
\begin{equation}\label{eq:back-stable-error}
	\frac{\lvert \tilde f(x) - f(x) \rvert}{\lvert f(x) \rvert} < c \kappa_{f}(x) \emachine.
\end{equation}
Furthermore, this bound composes when the underlying functions are well-conditioned: if $\tilde f$ and $\tilde g$ are both backward stable algorithms, then
\begin{equation*}
	\frac{\lvert \tilde f(\tilde g(x)) - f(g(x)) \rvert}{\lvert f(g(x)) \rvert} < c \kappa_{f}(g(x)) \kappa_{g}(x) \emachine.
\end{equation*}
Note that if $\kappa_{f \circ g}(x) \approx \kappa_{f}(g(x)) \kappa_{g}(x)$, this is the bound we would get for $\tilde f \circ \tilde g$ as a backward stable algorithm.
As a corollary, any calculation constructed from backward stable parts (such as elementary arithmetic) that exhibits large errors must have ill-conditioned steps.
An intuitive and quantifiable strategy for designing stable algorithms is to ensure that every step is as well-conditioned as possible.

We now turn our attention to representative functions in solid mechanics.
Most strain-energy functions and corresponding stress models in hyperelasticity contain expressions like $\log (1+x)$ or $\exp(x) - 1$, which are numerically unstable expressions when $x \approx 0$.
Numerical analysts proposed \cite{beebe2017mathematical, beebe2002computation}
\begin{align*}
	\logp(x) &= \log(1+x) & \expm(x) &= \exp(x) - 1,
\end{align*}
which gives full accuracy for small values of $x$.
These functions are specified by C99 and IEEE 754-2008 and widely available in core math libraries; see also \autoref{appendix:log1p}. \autoref{fig:log1p(x)-1-plus-x-m1} shows the $x, (1+x)-1$ and $\log(1+x), \logp(x)$ functions around 0. In both cases, large relative error is incurred by an algorithmic step that maps values near 1 to values near 0 via a function of derivative about 1, leading to unbounded condition number \eqref{eq:def-cond-number}.
Note that $\log(1+x)$ and $\exp(x) - 1$ are both well conditioned, but their direct evaluation is numerically unstable due to an ill-conditioned step.
We surveyed many open source finite element analysis packages and found that all contain numerically unstable formulations at some point due to phenomena explained above.

\subsection{Finite-strain mechanics}
Let $\bm X$ be the reference configuration and $\bm x = \bm X + \bm u$ be the current configuration expressed in terms of the displacement $\bm u$.
Wriggers \cite{wriggers2008nonlinear} discusses the displacement gradient
$$\bm H = \frac{\partial \bm u}{\partial \bm X}$$
and mentions that the Green-Lagrange strain can be expressed as
\begin{equation}\label{eq:green-lagrange-stable}
	\bE(\bH) = \frac 1 2 (\bH + \bH^T + \bH^{T} \bH)
\end{equation}
``in analytical investigations'', but notes ``This is actually not necessary when a numerical approach is applied.''
The usual presentation defines the deformation gradient $\bF = \bI + \bH$ and proceeds to the right Cauchy-Green tensors $\bC = \bF^{T} \bF$ and $\bE = \frac{1}{2} (\bC - \bI)$.
At small strain, $\bE$ is small despite $\bF$ and $\bC$ being of order 1, leading to instability as in the examples above.
For stable compressible formulations, one must also formulate expressions involving $J = \Det\bF$ and the first and second invariants of $\bC$ in a stable way.

Developers of FEAP \cite{feap2020} observed numerical stability issues when applying the standard approach at small strains and have added the ability to define user materials in terms of the displacement gradient $\bH$, which permits stable calculation of a strain-like quantity $\bb-\bI = \bH + \bH^{T} + \bH \bH^{T}$ \cite{govindjee2024personal,feap2020}.
That work, led by Bob Taylor circa 2011 and unpublished, reformulated several hyperelastic constitutive models in current configuration to use $J-1$ and $\bb - \bI$, with a Taylor series approximation of the volumetric terms when $\lvert J - 1 \rvert < 0.001$ \cite{govindjee2024personal}.
Other than FEAP, we are not aware of any software using stable formulations, nor of any stable formulations for initial configuration or any publications addressing these issues.
Moreover, user material interfaces in popular packages prevent stable formulations, e.g., by providing access to quantities such as $\bF$, $J$, and/or a Cauchy-Green tensor or its invariants instead of $\bH$, $J-1$, and/or nonlinear strain tensors.
This is indeed the case for Abaqus \texttt{UMAT} \cite{Abaqus2009}, which provides the deformation gradient $\bF$ and a strain increment, but does not provide direct access to the displacement gradient $\bH$ or nonlinear strain tensor, therefore numerical cancellation is inevitable when small strains appear within large strain formulations.
The Abaqus \texttt{UHYPER} interface is similarly unstable due to passing $J$ and invariants of the (isochoric) Cauchy-Green tensor.
FEBio \cite{maas2012febio}, a nonlinear finite element package for biomechanical applications, uses the deformation gradient $\bF$ even for defining the linear strain tensor.
Moreover, they offer a variety of hyperelasticity models that contains the subtraction, leading to loss of significance in small deformation regimes.
MoFEM \cite{mofemJoss2020}, an open source library for solving complex physics problems, uses the standard formulation given in continuum mechanics textbooks to describe the stress-strain relation of material.
This constitutive formulation contains $\log J$, leading to instability when $J \approx 1$.
Similarly, the Multiphysics Object-Oriented Simulation Environment (MOOSE) \cite{lindsay2022moose} defines, for example, the Neo-Hookean model using the standard (numerically unstable) formulation.
Albany-LCM  \cite{salinger2016albany}, a finite element code for analysis of multiphysics problems on unstructured meshes, expresses the hyperelastic models in terms of $J^2-1$, leading to catastrophic cancellation when $J \approx 1$.
\autoref{table:stable-constitutive} summarizes the stability properties of formulations used in well-known textbooks and production software.

\begin{table}[t]
	\caption{A stable constitutive implementation requires stable computation of strain $\bE$, $J = \Det\bF$, $\log J$, and shear expressions involving first and second invariants of $\bC$. Unstable formulations are widespread in books and software, with Ratel having the only stable formulation we are aware of.}
	\centering
	\begin{tabular}{l*{5}{c}}
		\toprule
		\multicolumn{1}{c}{Software/book} & \multicolumn{1}{c}{stable strain $\bE$} & \multicolumn{1}{c}{stable $J$} & \multicolumn{1}{c}{stable $\log J$} & \multicolumn{1}{c}{stable constitutive equation} \\
		\midrule
		FEAP \cite{feap2020} & \cmark  & \cmark & \cmark  & \cmark   \\
		FEBio \cite{maas2012febio} & \xmark  & \xmark & \xmark  & \xmark   \\
		Abaqus-\texttt{UMAT} \cite{Abaqus2009} & \xmark  & \xmark & \xmark  & \xmark   \\
		MOOSE \cite{lindsay2022moose} & \xmark  & \xmark & \xmark  & \xmark   \\
		Albany-LCM \cite{salinger2016albany} & \xmark  & \xmark & \xmark  & \xmark   \\
		LifeV \cite{bertagna2017lifev} & \xmark & \xmark & \xmark & \xmark \\
		MoFEM \cite{mofemJoss2020} & \xmark & \xmark & \xmark & \xmark \\
		Ratel \cite{ratelusermanual} & \cmark  & \cmark & \cmark  & \cmark   \\
		Holzapfel \cite{holzapfel2000nonlinear}  & \xmark & \xmark & \xmark  & \xmark   \\
		Wriggers \cite{wriggers2008nonlinear}  & \xmark  & \xmark & \xmark & \xmark \\
		\bottomrule
	\end{tabular}
	\label{table:stable-constitutive}
\end{table}

The paper proceeds as follows: \autoref{constitutive} develops stable formulations for common hyperelastic models in coupled and decoupled (isochoric and volumetric split), \autoref{ad} demonstrates that stable stress expressions can be derived using algorithmic differentiation (AD) so long as care is taken in the strain energy formulation, \autoref{axial} gives an illustrative numerical test, and \autoref{conclusion} concludes with outlook toward inelastic models. Details of the numerical procedure for evaluating stability are given in \autoref{num-eval-stability}.
All figures exhibited here are created using the open source Julia programming language. Comprehensive numerical experiments and figures are provided in the executable supplement \cite{shakeri2024stabledemo} for those readers who are interested in exploring of all given formulations here.
While this study presents some of the most well-known hyperelastic models and their stable formulation, the approach is general and we provide guidance for applying it to other material models.

\section{Constitutive equations}\label{constitutive}
The constitutive behavior for hyperelastic materials is characterized by a strain energy density function $\psi$.
For isotropic materials in initial configuration, $\psi$ is typically defined by either the principal invariants $\{\mathbb{I}_1, \mathbb{I}_2, \mathbb{I}_3 \}$ of right Cauchy-Green tensor $\bC = \bF^T \bF$ or the principal stretches $\{\lambda_1, \lambda_2, \lambda_3 \}$ of the SPD matrix $\bU$ where $\bR\bU = \bF$ is the polar decomposition.
In the following we discuss the most common coupled and decoupled representation of the strain energies and the associated constitutive equations that are employed frequently in the literature.

\subsection{Coupled strain energy}
Coupled strain energy functionals $\psi$ are written in terms of invariants without an isochoric-volumetric split.
In the linear regime, this corresponds to use of shear modulus $\mu$ and first Lamé parameter $\firstlame$, with the standard (not deviatoric) infinitesimal strain tensor $\bm\varepsilon$.
For the general form of coupled strain energy
\begin{equation}\label{general-strain-energy}
	\psi \left( \bm{C} \right) = \psi \left( \mathbb{I}_1, \mathbb{I}_2, \mathbb{I}_3 \right),
\end{equation}
with the invariants
\begin{equation}\label{invariants}
	\begin{split}
		\mathbb{I}_1 (\bm{C}) &= \trace \bm{C}, \\
		\mathbb{I}_2 (\bm{C}) &= \frac 1 2 \left( \mathbb{I}_1^2 - \bm{C} \tcolon \bm{C} \right) \\
		\mathbb{I}_3 (\bm{C}) &= \Det{\bm{C}}.
	\end{split}
\end{equation}
we can determine the constitutive equations by taking the gradient of strain energy as
\begin{equation}\label{strain-energy-grad}
	\bm{S} = \frac{\partial \psi}{\partial \bm{E}}= 2 \frac{\partial \psi}{\partial \bm{C}} = 2 \sum_{i=1}^3 \frac{\partial \psi}{\partial \mathbb{I}_i} \frac{\partial \mathbb{I}_i}{\partial \bm{C}},
\end{equation}
where $\bm{S}$ is the general form of the stress relation in initial configuration and
\begin{equation}\label{invariants-derivative}
	\begin{split}
		\frac{\partial \mathbb{I}_1}{\partial \bm{C}} &= \bm{I}, \\
		\frac{\partial \mathbb{I}_2}{\partial \bm{C}} &=  \mathbb{I}_1 \bm{I} - \bm{C}\\
		\frac{\partial \mathbb{I}_3}{\partial \bm{C}} &= \mathbb{I}_3 \bm{C}^{-1}.
	\end{split}
\end{equation}
In the following we introduce the stable formulation for two well-known hyperelastic constitutive equations.

\subsubsection{Neo-Hookean model}\label{coupled-nh}
One of the simplest hyperelastic models is the Neo-Hookean model, given by
\begin{equation} \label{NH-energy}
	\psi(\bC) = \frac{\firstlame}{4} \left( J^2 - 1 -2 \log J \right) -\mu \log J +\frac{\mu}{2} \left(\mathbb{I}_1 - 3 \right)
\end{equation}
where $ J =\Det\bF = \sqrt{\Det\bC}$. The first term is a convex choice satisfying limit conditions \cite{wriggers2008nonlinear,Doll2000OnTD} while the second is a structural necessity for coupled strain energy formulations. The second Piola-Kirchhoff stress is derived according to \eqref{invariants-derivative}
\begin{equation}\label{NH-S-unstable}
	\bS = 2\frac{\partial \psi}{\partial \bC} = \frac{\partial \psi}{\partial \bE} = \frac{\firstlame}{2} \left( J^2 - 1 \right)\bC^{-1} + \mu \left(\bI - \bC^{-1} \right),
\end{equation}
where $\bC = \bI + 2 \bE$. Both terms in \eqref{NH-S-unstable} are numerically unstable at small strain $\bC \approx I$. The first term is also unstable at large strain when $J \approx 1$, which is typical in finite-strain analysis of nearly-incompressible materials. Similar to $\logp$, we need a formulation that avoids direct computation of $J$ in $J^{2} - 1 = (J-1)(J+1)$.

Consider the 2-dimensional case
$$
\bF = \bI + \bH = \bI + \nabla_X \bm{u} = \begin{bmatrix}
	1 + u_{1,1} & u_{1,2} \\
	u_{2,1} & 1 + u_{2,2}
\end{bmatrix},
$$
and let $\Jm = J-1$ be computed by the stable expression
\begin{equation}\label{Jm1}
	\mathtt{J_{-1}} = u_{1,1} + u_{2,2} + u_{1,1} u_{2,2} - u_{1,2} u_{2,1}.
\end{equation}

\begin{remark}
	Ideally, one would like a proof of backward stability for every constitutive model, but doing so is tedious and we believe offers no great insight. In light of \eqref{eq:back-stable-error}, it is revealing to plot the forward relative error and observe that stable algorithms provide errors that are uniformly of order $\emachine$ (sometimes written ``in $O(\emachine)$'') independent of the input argument $\bH$ (or a strain).
	There is a caveat: in the case of tiny perturbations of a pure rotation $\bH \approx \bm Q - \bI$ for an orthogonal matrix $\bm Q$, the strain is nearly zero and thus stress will also be nearly zero.
	For example, the Green-Lagrange strain $\bE(\bH) = \frac 1 2 (\bF^T \bF - \bI) \approx \frac 1 2 (\bm Q^{T} \bm Q - \bI) = \bm 0$ is a function with unbounded condition number, thus even a backward stable algorithm such as \eqref{eq:green-lagrange-stable} will have unbounded error in $\bE$.
	One might consider high/mixed precision or a decomposition of the state variable $\bm u$ into a finite rotation and (possibly small) perturbation \cite{casey1992infinitesimal} so that strain can be made a well-conditioned function of the state for problems that require small-strain stability for large motions of floating bodies.
	While we consider $\bH$ as the input in our numerical experiments (to be agnostic over strain measures), we will not construct cases of nearly-pure rotations.
\end{remark}

\begin{remark}
	Apart from the incompressible limit, hyperelastic constitutive models are well-conditioned, thus \eqref{eq:back-stable-error} ensures the error will be in $O(\emachine)$.
	Although various displacement-only formulations are common in engineering practice, mixed methods are necessary for well-conditioned finite element formulations in the incompressible limit.
\end{remark}

\begin{figure}
	\centering
	\includegraphics[width=.7\linewidth]{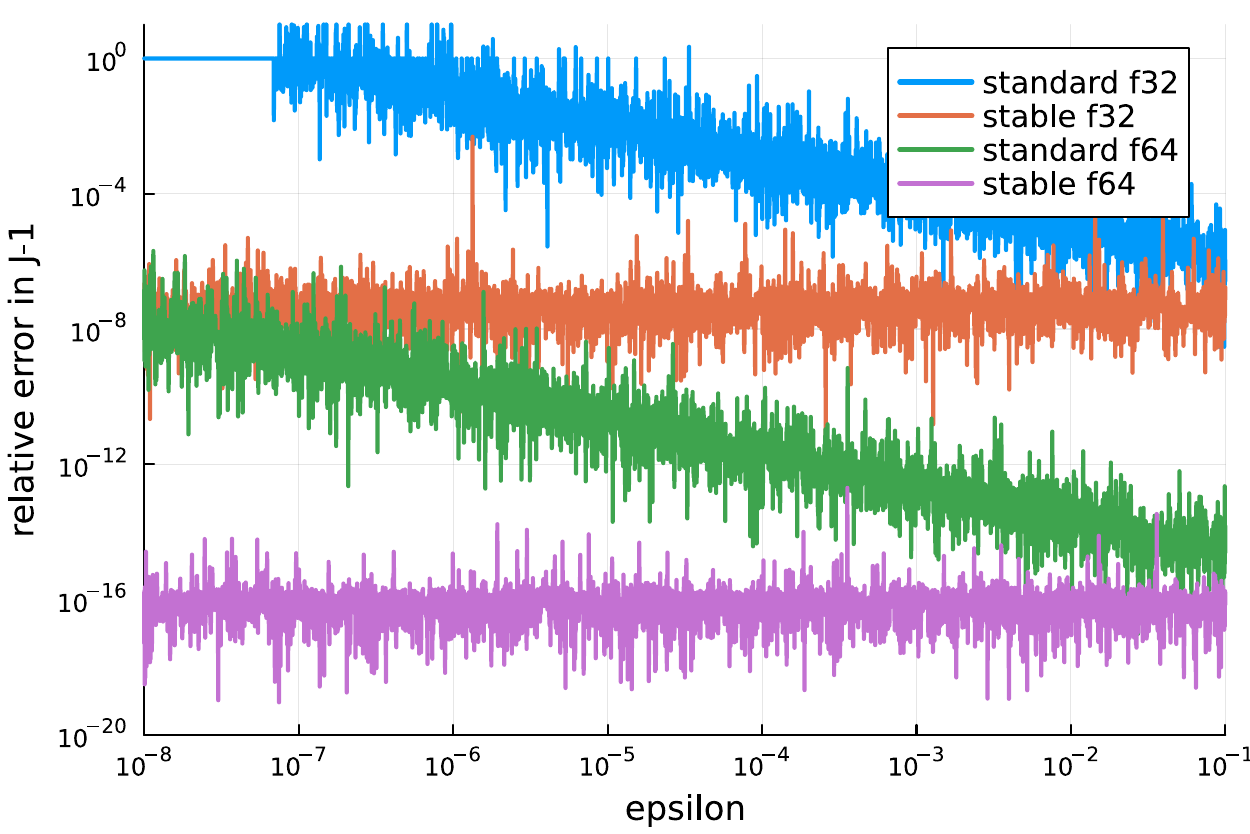}
	\caption{Relative error of standard computation of $J-1$ and its stable way $\mathtt{J_{-1}}$.} \label{fig1:Jm1}
\end{figure}

Using a 3-dimensional analog of \eqref{Jm1} described in \autoref{num-eval-stability} and numerical procedure explained therein, \autoref{fig1:Jm1} shows that the relative error in $\Jm$ is $O(\emachine)$ independent of the magnitude of the displacement gradient $\bH$ and strain, while $J-1$ loses digits of accuracy at small strain.
For instance, for strain at order $10^{-8}$ we cannot trust any digits in the final computed $J-1$ by standard approach with single precision, and only half of the digits in double precision are accurate.
Replacing $J^{2}-1$ with $\Jm \left(\Jm + 2 \right)$ is sufficient to stabilize the first term of \eqref{NH-S-unstable} and get all 8 and 16 digits accurate with single and double precision for all range of elasticity strain. 

In some constitutive models, $\log J$ appears directly in the expression for $\bS$, and also has a huge condition number $\kappa_{\log}(J) = 1 / \log J $ when $ J \approx 1$.
In such cases, one naturally achieves stability via
\begin{equation}\label{logJ-stable}
	\log J = \logp \left(\mathtt{J_{-1}} \right).
\end{equation}

The Neo-Hookean constitutive equation \eqref{NH-S-unstable} is also unstable due to cancellation in the second term when $\bC^{-1} \approx \bI$, and it can be fixed by replacing $\bI - \bC^{-1} = \bC^{-1}(\bC - \bI) = 2 \bC^{-1} \bE$. Thus, an equivalent stable form of \eqref{NH-S-unstable} is
\begin{equation}\label{NH-S-stable}
	\bS = \frac{\firstlame}{2} \mathtt{J_{-1}} \left(\mathtt{J_{-1}} + 2\right)\bC^{-1} + 2 \mu \bC^{-1} \bE .
\end{equation}
The relative error in $\bS$ using the standard \eqref{NH-S-unstable} and stable \eqref{NH-S-stable} expressions yields a figure indistinguishable from \autoref{fig1:Jm1}.
Although we elide figures for each constitutive model and stress, the supplement \cite{shakeri2024stabledemo} contains numerical demonstration that every expression we call unstable or stable exhibits errors indistinguishable from the corresponding line in \autoref{fig1:Jm1}.

The Neo-Hookean constitutive model can also be evaluated in current configuration, classically written in terms of $J$ and the left Cauchy-Green tensor $\bb = \bF \bF^T$.
One can derive an expression for a current configuration stress (we use Kirchhoff stress $\btau$) starting with a strain energy density $\psi(\bb)$ or (equivalently) pushing $\bS$ of \eqref{NH-S-unstable} forward via
\begin{equation}\label{NH-tau-unstable}
	\btau = \bF\bS\bF^T = \frac{\firstlame}{2} \left( J^2 - 1 \right)\bI + \mu \left(\bb - \bI \right).
\end{equation}
A stable form analogous to \eqref{NH-S-stable} is
\begin{equation}\label{NH-tau-stable}
	\btau = \frac{\firstlame}{2} \mathtt{J_{-1}} \left(\mathtt{J_{-1}} + 2\right) \bI + 2 \mu \be,
\end{equation}
where we call $\be = \frac 1 2 (\bb - \bI) = \frac 1 2 (\bH + \bH^T + \bH \bH^T)$ the Green-Euler strain tensor.
(It is the $n=2$ current configuration Seth-Hill strain measure \cite{holzapfel2000nonlinear}.)

\subsubsection{Mooney-Rivlin model}
Compared to the Neo-Hookean model \eqref{NH-energy}, the coupled energy function for Mooney-Rivlin model depends on the additional invariant $\mathbb{I}_2(\bC) = \frac 1 2 \left(\mathbb{I}_1^2(\bC) - \bC \tcolon \bC \right)$ where $\bC \tcolon \bC = \trace(\bC^{T} \bC)$ is the Frobenius inner product. The coupled Mooney-Rivlin strain energy is given by \cite{holzapfel2000nonlinear}
\begin{equation}\label{MR-energy}
	\psi(J, \bC) = \frac{\firstlame}{4} \left( J^2 - 1 -2 \log J \right) -\left( \mu_1 + 2\mu_2 \right) \log J +\frac{\mu_1}{2} \left( \mathbb{I}_1(\bC) - 3 \right) + \frac{\mu_2}{2} \left( \mathbb{I}_2(\bC) - 3 \right),
\end{equation}
where $\mu_{1}$ and $\mu_{2}$ are parameters that must be experimentally determined and $\mu = \mu_{1}+\mu_{2}$ is the shear modulus in the linear regime.
The second Piola-Kirchhoff tensor is
\begin{equation}\label{MR-S-unstable}
	\bS = \frac{\firstlame}{2} \left( J^2 - 1 \right)\bC^{-1} + \mu_1 \left(\bI - \bC^{-1} \right) + \mu_2 \left(\mathbb{I}_1(\bC)\bI - 2\bC^{-1} - \bC \right).
\end{equation}
Similar to the Neo-Hookean model, we write the stable form in terms of Green-Lagrange strain,
\begin{equation}\label{MR-S-stable}
	\bS = \frac{\firstlame}{2} \mathtt{J_{-1}} \left(\mathtt{J_{-1}} + 2\right)\bC^{-1} + 2 \left( \mu_1 + 2 \mu_2 \right)\bC^{-1} \bE + 2 \mu_2 \left(\mathbb{I}_1(\bE)\bI - \bE \right).
\end{equation}
In addition, the Kirchhoff stress tensor for Mooney-Rivlin in current configuration is classically written
\begin{equation}\label{MR-tau-unstable}
	\btau = \frac{\firstlame}{2} \left( J^2 - 1 \right)\bI + \mu_1 \left(\bb - \bI \right) + \mu_2 \left(\mathbb{I}_1(\bb) \bb - 2\bI - \bb^2 \right)
\end{equation}
and its stable form in terms of Green-Euler strain is
\begin{equation}\label{MR-tau-stable}
	\btau = \frac{\firstlame}{2} \mathtt{J_{-1}} \left( \mathtt{J_{-1}} + 2\right)\bI +2 \left( \mu_1 + 2 \mu_2 \right)\be + 2 \mu_2 \left(\mathbb{I}_1(\be)\bI - \be \right) \bb.
\end{equation}

\subsection{Decoupled strain energy}
Decoupled strain energy formulations decompose the strain into volumetric and deviatoric parts. In the linear regime, such formulations use the shear modulus $\mu$ and bulk modulus $\bulk$, with the deviatoric strain $\bm\varepsilon' = \bm\varepsilon - \frac 1 3 (\trace \bm \varepsilon) \bI$.
In the case of rubber-like materials, the bulk modulus $\bulk$ is orders of magnitude larger than the shear modulus, $\bulk \gg \mu$.
For finite strain, the strain energy function is split into volumetric and isochoric (volume-preserving) parts.
The first step is multiplicative decomposition of the deformation gradient as
\begin{equation}\label{modified-deformation-gradient}
	\bF = (J^{1/3} \bI) \bar{\bF} = J^{1/3} \bar{\bF}
\end{equation}
where $J^{1/3} \bI$ describes the purely volumetric deformation and $\bar{\bF}$ captures the isochoric  since $\Det{\bar{\bF}}  = \Det{J^{-1/3} \bF}  = 1$ .
Similarly we can decompose right Cauchy-Green tensor
\begin{equation}\label{modified-right-cauchy-green}
	\bar{\bC} = \bar{\bF}^T \bar{\bF} = J^{-2/3} \bC
\end{equation}
where $\bar{\bm{F}}$ and $\bar{\bm{C}}$ are known as the modified deformation gradient and modified right Cauchy-Green tensor, respectively. Moreover, the modified principal stretches $\bar{\lambda}_i$ and modified invariants $\mathbb{\bar{I}}_i$ can be defined
\begin{equation}\label{modified-stretch}
	\bar{\lambda}_i = J^{-1/3} \lambda_i, \, i=1,2,3
\end{equation}
\begin{align}\label{modified-invariants}
	\mathbb{\bar{I}}_1 &= J^{-2/3}\mathbb{I}_1, & \mathbb{\bar{I}}_2 &= J^{-4/3}\mathbb{I}_2, & \mathbb{\bar{I}}_3 &= \left(J^{-2/3}\right)^3 \mathbb{I}_3 = 1.
\end{align}
The general decoupled strain energy is a sum of volumetric and isochoric parts
\begin{equation}\label{general-strain-energy-decoupled}
	\psi \left( J, \bar{\bm{C}} \right) = \psi_{\text{vol}}\left(J \right) + \psi_{\text{iso}} \left(\bar{\bC} \right),
\end{equation}
leading to an additive decomposition of the second Piola-Kirchhoff stress
\begin{equation}\label{S-iso-vol}
	\bS =\frac{\partial \psi}{\partial \bE} = \frac{\partial \psi_{\text{vol}}}{\partial J} \frac{\partial J}{\partial \bE} + \frac{\partial \psi_{\text{iso}}}{\partial \bE} = \bS_{\text{vol}} + \bS_{\text{iso}}.
\end{equation}
When the isochoric strain energy is written in terms of modified invariants $\psi_{\text{iso}}(\mathbb{\bar I}_1, \mathbb{\bar I}_2)$ (note that $\mathbb{\bar I}_{3}=1$ uniformly), the isochoric stress satisfies
\begin{equation}\label{S-iso-general}
	\bS_{\text{iso}} = \frac{\partial \psi_{\text{iso}}}{\partial \bE} = 2 \frac{\partial \psi_{\text{iso}}}{\partial \bm{C}} = 2 \sum_{i=1}^3 \frac{\partial \psi_{\text{iso}}}{\partial \mathbb{\bar{I}}_i} \frac{\partial \mathbb{\bar{I}}_i}{\partial \bm{C}}
\end{equation}
with
\begin{equation} \label{modified-invariants-derivative}
	\begin{split}
		\frac{\partial \mathbb{\bar{I}}_1}{\partial \bm{C}} &= \frac{\partial (J^{-2/3}\mathbb{I}_1)}{\partial \bm{C}} = J^{-2/3} \left( \bm{I} -\frac{1}{3} \mathbb{I}_1 \bm{C}^{-1}\right) \\
		\frac{\partial \mathbb{\bar{I}}_2}{\partial \bm{C}} &=  \frac{\partial (J^{-4/3}\mathbb{I}_2)}{\partial \bm{C}} = J^{-4/3} \left( \mathbb{I}_1 \bm{I} - \bm{C} - \frac{2}{3} \mathbb{I}_2 \bm{C}^{-1}\right), \\
	\end{split}
\end{equation}
where we have used $\frac{\partial J}{\partial \bm{C}} = \frac{1}{2} J \bm{C}^{-1}$. As discussed for \eqref{NH-S-unstable}, there are many empirical forms for the pure volumetric part (or bulk term)
of the strain-energy, which should be convex and satisfy physical constrains \cite{Doll2000OnTD, MOERMAN2020474}. We consider one such form,
\begin{equation}\label{energy-vol}
	\psi_{\text{vol}} = \frac{\bulk}{4} \left( J^2 - 1 -2 \log J \right),
\end{equation}
but similar principles will give stable formulations for others. The volumetric stress can be defined as
\begin{equation}\label{S-vol-general}
	\bS_{\text{vol}} = -p \frac{\partial J}{\partial \bE} = -p J \bC^{-1}
\end{equation}
where we have used the definition of hydrostatic pressure $p = - \frac{\partial \psi_{\text{vol}}}{\partial J}$, which for \eqref{energy-vol} can be stably computed via
\begin{equation}\label{pressure-stable}
	p = -\frac{\bulk}{2 J} \left(J^2 - 1 \right) = -\frac{\bulk}{2 J} \mathtt{J_{-1}} \left(\mathtt{J_{-1}} + 2 \right).
\end{equation}
While $p$ depends on the form $\psi_{\text{vol}}(J)$, the volumetric stress \eqref{S-vol-general} is numerically stable and always the same expression in terms of $p$. In the decoupled framework, the only salient difference between the Neo-Hookean, Mooney-Rivlin, and Ogden models is in their isochoric part, thus we focus on stable expressions for isochoric stresses $\bS_{\text{iso}}$ and $\btau_{\text{iso}}$.

\subsubsection{Neo-Hookean model}
The decoupled strain energy density for the Neo-Hookean model is
\begin{equation}\label{NH-energy-decoupled}
	\psi(J, \bar{\bC}) = \psi_{\text{vol}}(J) + \frac{\mu}{2} \left( \bar{\mathbb{I}}_1 - 3 \right)
\end{equation}
The isochoric second Piola-Kirchhoff stress is a straightforward application of \eqref{modified-invariants-derivative},
\begin{equation}\label{NH-S-unstable-iso}
	\bS_{\text{iso}} = \mu J^{-2/3}\left(\bm{I} - \frac{1}{3} \mathbb{I}_1 \bm{C}^{-1} \right).
\end{equation}
Using the relation $\mathbb{I}_1(\bC) = 3 + 2 \mathbb{I}_1(\bE)$, the numerically stable form of \eqref{NH-S-unstable-iso} can be written as
\begin{equation}\label{NH-S-stable-iso}
	\bS_{\text{iso}} = 2 \mu J^{-2/3} \bC^{-1} \left(\bE -\frac{1}{3}\mathbb{I}_1(\bE) \bI \right) = 2 \mu J^{-2/3} \bC^{-1} \bE_{\text{dev}},
\end{equation}
which makes use of the deviatoric Green-Lagrange strain $\bE_{\text{dev}} = \bE -\frac{1}{3}\mathbb{I}_1(\bE) \bI$.

In current configuration, the isochoric Kirchhoff stress is
\begin{equation}\label{NH-tau-unstable-iso}
	\btau_{\text{iso}} = \bF \bS_{\text{iso}} \bF^T = \mu J^{-2/3}\left(\bb - \frac{1}{3}\mathbb{I}_1(\bb) \bm{I} \right)
\end{equation}
and its equivalent stable form is
\begin{equation}\label{NH-tau-stable-iso}
	\btau_{\text{iso}} = 2\mu J^{-2/3}\left(\be - \frac{1}{3}\mathbb{I}_1(\be) \bI \right) = 2\mu J^{-2/3} \be_{\text{dev}},
\end{equation}
where we have used $\mathbb{I}_1(\bb) = 3 + 2\mathbb{I}_1(\be)$.

\subsubsection{Mooney-Rivlin model}
For the Mooney-Rivlin model, decoupled strain energy density is given by
\begin{equation}\label{MR-energy-decoupled}
	\psi(J, \bar{\bC}) = \psi_{\text{vol}}(J) +\frac{\mu_1}{2} \left( \bar{\mathbb{I}}_1 - 3 \right) + \frac{\mu_2}{2} \left( \bar{\mathbb{I}}_2 - 3 \right)
\end{equation}
The isochoric second Piola-Kirchhoff stress can be written as \eqref{modified-invariants-derivative}
\begin{equation}\label{MR-S-unstable-iso}
	\bS_{\text{iso}} =\mu_1 J^{-2/3}\left(\bI - \frac{1}{3} \mathbb{I}_1(\bC) \bC^{-1} \right) + \mu_2 J^{-4/3}\left(\mathbb{I}_1(\bC)\bm{I} - \bC - \frac{2}{3} \mathbb{I}_2(\bC) \bC^{-1} \right)
\end{equation}
and its stable form is
\begin{equation}\label{MR-S-stable-iso}
	\begin{split}
		\bS_{\text{iso}} &= 2 \mu_1 J^{-2/3} \bC^{-1} \left(\bE -\frac{1}{3}\mathbb{I}_1 (\bE) \bI \right) + 2 \mu_2 J^{-4/3} \left(\mathbb{I}_1(\bE) \bI - \bE \right) + 4 \mu_2 J^{-4/3} \bC^{-1} \left(\bE - \frac{2}{3} \mathbb{I}_1(\bE) \bI - \frac{2}{3} \mathbb{I}_2(\bE) \bI \right) \\
		&= 2 \left( \mu_1 J^{-2/3} + 2 \mu_2 J^{-4/3} \right) \bC^{-1} \bE_{\text{dev}} + 2 \mu_2 J^{-4/3} \left(\mathbb{I}_1(\bE) \bI - \bE \right) - \frac{4}{3} \mu_2 J^{-4/3} \left( \mathbb{I}_1(\bE) + 2\mathbb{I}_2(\bE) \right) \bC^{-1},
	\end{split}
\end{equation}
where we have used $\mathbb{I}_1(\bC) = 3 + 2\mathbb{I}_1(\bE), \, \mathbb{I}_2(\bC) = 3 + 4\mathbb{I}_1(\bE) + 4\mathbb{I}_2(\bE)$.

In current configuration we have
\begin{equation}\label{MR-tau-unstable-iso}
	\btau_{\text{iso}} = \mu_1 J^{-2/3}\left(\bb - \frac{1}{3}\mathbb{I}_1(\bb) \bI \right) + \mu_2 J^{-4/3}\left(\mathbb{I}_1(\bb) \bb - \bb^2 - \frac{2}{3}\mathbb{I}_2(\bb) \bI \right)
\end{equation}
and the stable form of \eqref{MR-tau-unstable-iso} is
\begin{equation}\label{MR-tau-stable-iso}
	\begin{split}
		\btau_{\text{iso}} &= 2\mu_1 J^{-2/3}\left(\be - \frac{1}{3}\mathbb{I}_1 (\be) \bI \right) + 2 \mu_2 J^{-4/3} \left(\mathbb{I}_1(\be) \bI - \be \right) \bb + 4 \mu_2 J^{-4/3} \left(\be - \frac{2}{3} \mathbb{I}_1(\be) \bI - \frac{2}{3} \mathbb{I}_2(\be) \bI \right) \\
		&= 2 \left( \mu_1 J^{-2/3} + 2 \mu_2 J^{-4/3} \right) \be_{\text{dev}} + 2 \mu_2 J^{-4/3} \left(\mathbb{I}_1(\be) \bI - \be \right) \bb - \frac{4}{3} \mu_2 J^{-4/3} \left( \mathbb{I}_1(\be) + 2\mathbb{I}_2(\be) \right) \bI
	\end{split}
\end{equation}
where we have used $\mathbb{I}_1(\bb) = 3 + 2\mathbb{I}_1(\be), \, \mathbb{I}_2(\bb) = 3 + 4\mathbb{I}_1(\be) + 4\mathbb{I}_2(\be)$.

\subsubsection{Ogden model}
The postulated decoupled strain energy for Ogden model is a function of the modified principal stretches \eqref{modified-stretch} as \cite{wriggers2008nonlinear, holzapfel2000nonlinear}

\begin{equation}\label{Og-energy-decoupled}
	\psi(J, \bar{\lambda}_i) = \psi_{\text{vol}}(J) + \psi_{\text{iso}}(\bar{\lambda}_i) = \psi_{\text{vol}}(J) + \sum_{i=1}^3\bar{\omega}(\bar{\lambda}_i) \quad \text{with} \quad
	\bar{\omega}(\bar{\lambda}_i) = \sum_{j=1}^{N}\frac{\mu_j}{\alpha_j}\left(\bar{\lambda}_i^{\alpha_j} - 1 \right)
\end{equation}
where the  parameters $\mu_j$ and $\alpha_j$ have to be determined from experiments. 
In the linearized regime, all hyperelastic models reduce to linear elasticity, where the shear modulus $\mu$ satisfies
\begin{equation}
	2\mu = \sum_{j=1}^N \mu_j \alpha_j \quad \text{with} \quad \mu_j \alpha_j > 0
\end{equation}
in terms of the Ogden parameters.
The isochoric second Piola-Kirchhoff stress can be written as
\begin{equation}\label{Og-S-unstable-iso}
	\bS_{\text{iso}}
	=2\frac{\partial \psi_{\text{iso}} (\bar{\lambda}_i)}{\partial \bC}
	=\sum_{i=1}^{3} 2 \frac{\partial \psi_{\text{iso}}}{\partial \lambda_i^2}
	\frac{\partial \lambda_i^2}{\partial \bC}
	=\sum_{i=1}^{3} \frac{1}{\lambda_i}\frac{\partial \psi_{\text{iso}}}{\partial \lambda_i} \bN_i \bN_i^T = \sum_{i=1}^{3} s_i \bN_i \bN_i^T
\end{equation}
where an eigenvector $\bN_i$ is computed by $\bC\bN_i = \lambda_i^2 \bN_i$. By employing $\frac{\partial J}{\partial {\lambda}_i} = J {\lambda}_i^{-1}$ we have
\begin{equation}\label{si-unstable}
	s_i = \frac{1}{\lambda_i}\frac{\partial \psi_{\text{iso}}}{\partial \lambda_i}
	=\frac{1}{\lambda_i}\frac{\partial \bar{\lambda}_k}{\partial \lambda_i}
	\frac{\partial \psi_{\text{iso}}}{\partial \bar{\lambda}_k}
	=\frac{J^{-1/3}}{\lambda_i} \left(\delta_{ik} -\frac{1}{3}\lambda_i^{-1}\lambda_k \right) \frac{\partial \psi_{\text{iso}}}{\partial \bar{\lambda}_k}
	=\frac{J^{-1/3}}{\lambda_i} \left(\delta_{ik} -\frac{1}{3}\bar{\lambda}_i^{-1}\bar{\lambda}_k \right) \frac{\partial \psi_{\text{iso}}}{\partial \bar{\lambda}_k}
\end{equation}
in which
\begin{equation}
	\frac{\partial \psi_{\text{iso}}}{\partial \bar{\lambda}_k} = \frac{\partial \bar{\omega} (\bar{\lambda}_k)}{ \partial \bar{\lambda}_k} = \sum_{j=1}^N \mu_j \bar{\lambda}_k^{\alpha_j - 1}
\end{equation}

To derive an equivalent numerically stable form of \eqref{Og-S-unstable-iso} we rewrite the $\psi_{\text{iso}}$ by substituting $\bar{\lambda}_i = J^{-1/3}\lambda_i$ as
\begin{equation} \label{psi-iso}
	\psi_{\text{iso}}({\lambda}_i) = \sum_{j=1}^N \frac{\mu_j}{\alpha_j} \left[\left({\lambda}_1^{\alpha_j} + {\lambda}_2^{\alpha_j} + {\lambda}_3^{\alpha_j} \right)J^{-\alpha_j/3}  -3 \right]
\end{equation}
and compute its derivative as
\begin{equation}\label{s1-stable}
	\begin{split}
		s_1 = \frac{1}{\lambda_1}\frac{\partial \psi_{\text{iso}}}{\partial \lambda_1} &= \frac{1}{\lambda_1} \sum_{j=1}^N \frac{\mu_j}{\alpha_j} \left[ \left(\alpha_j \lambda_1^{\alpha_j -1} \right)J^{-\alpha_j/3} - \frac{\alpha_j}{3} J^{-\alpha_j/3} \lambda_1^{-1} \left( {\lambda}_1^{\alpha_j} + {\lambda}_2^{\alpha_j} + {\lambda}_3^{\alpha_j} \right)\right] \\
		&= \frac{1}{\lambda_1} \sum_{j=1}^N \mu_j \left[ \lambda_1^{\alpha_j -1} - \frac{1}{3} \lambda_1^{-1} \left( {\lambda}_1^{\alpha_j} + {\lambda}_2^{\alpha_j} + {\lambda}_3^{\alpha_j} \right)\right] J^{-\alpha_j/3} \\
		&= \frac{1}{\lambda_1^2} \sum_{j=1}^N \frac{\mu_j}{3} \left[ 2\lambda_1^{\alpha_j} - {\lambda}_2^{\alpha_j} - {\lambda}_3^{\alpha_j} \right] J^{-\alpha_j/3} \\
		&= \frac{1}{\lambda_1^2} \sum_{j=1}^N \frac{\mu_j}{3} \left[ 2\left(\lambda_1^{\alpha_j} -1 \right) - \left({\lambda}_2^{\alpha_j} -1 \right) - \left({\lambda}_3^{\alpha_j} -1 \right)\right] J^{-\alpha_j/3} \\
		&=\frac{1}{\lambda_1^2} \sum_{j=1}^N \frac{\mu_j}{3} \left[ 2(e^{\alpha_j \ell_1}-1) - (e^{\alpha_j \ell_2}-1) - (e^{\alpha_j \ell_3}-1) \right] J^{-\alpha_j/3} \\
		&=\frac{1}{1 + 2\lambda_1^E} \sum_{j=1}^N \frac{\mu_j}{3} \left[ 2\operatorname{\tt expm1}(\alpha_j \ell_1) - \operatorname{\tt expm1}(\alpha_j \ell_2) - \operatorname{\tt expm1}(\alpha_j \ell_3) \right] J^{-\alpha_j/3}
	\end{split}
\end{equation}
where $\lambda_i^E$ is the eigenvalue of the strain tensor $\bE$ and $\ell_i = \log \lambda_i = \frac 1 2 \logp(2\lambda_i^E)$.
Note that $\bE$ has the same eigenvectors as $\bC$ without the loss of precision incurred by rounding in the formation of $\bC$.
Following the above approach, we have
\begin{equation}\label{s23-stable}
	\begin{split}
		s_2 = \frac{1}{\lambda_2}\frac{\partial \psi_{\text{iso}}}{\partial \lambda_2} &=
		\frac{1}{1 + 2\lambda_2^E} \sum_{j=1}^N \frac{\mu_j}{3} \left[ -\operatorname{\tt expm1}(\alpha_j \ell_1) + 2\operatorname{\tt expm1}(\alpha_j \ell_2) - \operatorname{\tt expm1}(\alpha_j \ell_3) \right] J^{-\alpha_j/3} \\
		s_3 = \frac{1}{\lambda_3}\frac{\partial \psi_{\text{iso}} }{\partial \lambda_3} &=
		\frac{1}{1 + 2\lambda_3^E} \sum_{j=1}^N \frac{\mu_j}{3} \left[ -\operatorname{\tt expm1}(\alpha_j \ell_1) - \operatorname{\tt expm1}(\alpha_j \ell_2) + 2\operatorname{\tt expm1}(\alpha_j \ell_3) \right] J^{-\alpha_j/3}.
	\end{split}
\end{equation}

\begin{figure}[t]
	\centering
	\includegraphics[width=.7\linewidth]{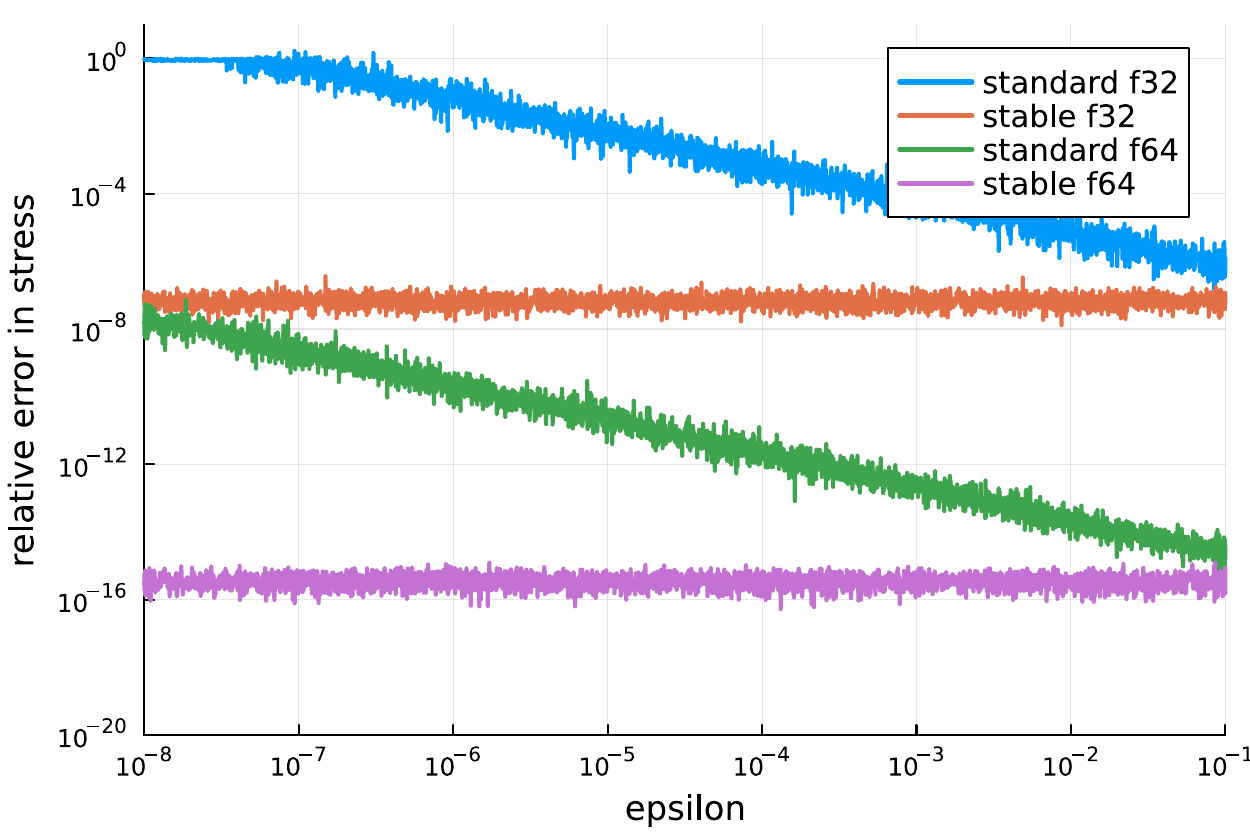}
	\caption{Relative error of standard computation of isochoric second Piola-Kirchhoff stress for Ogden model and its stable form} \label{fig2:Ogden}
\end{figure}

Substituting the new definition of $s_i$ \eqref{s1-stable} and \eqref{s23-stable} into \eqref{Og-S-unstable-iso} provides a stable form of the isochoric second Piola-Kirchhoff stress for the Ogden model.
Relative error for the standard \eqref{si-unstable} and stable \eqref{s1-stable}-\eqref{s23-stable} form of the isochoric second Piola-Kirchhoff stress \eqref{Og-S-unstable-iso} for the Ogden model is plotted in \autoref{fig2:Ogden}. The behavior for single and double precision is similar to \autoref{fig1:Jm1}.
The proposed stable Ogden formulation has relative error up to $\emachine$, while the standard approach loses $m$ digits for $\epsilon = 10^{-m}$.

Note that for a symmetric, real-valued $3 \times 3$ strain tensor $\bE$, the standard closed-form (and iterative) eigensolvers are susceptible to loss of accuracy in the eigenvectors, in the sense that $\bE = \bN \Lambda^{E} \bN^{T}$ may not be satisfied to a relative error of $\emachine$ in the presence of near- or exactly-repeated eigenvalues. To avoid this problem, we recommend computing the eigenvalues $\lambda_i^E$ and corresponding eigenvectors $\bN_{i}$ using the stable algorithm proposed by Harari and Albocher \cite{harari2023computation}, which has a relative accuracy of $\emachine$.

\subsection{Hencky strain}

Standard material logarithmic or Hencky strain is defined as

\begin{equation} \label{EH-unstable}
	\bE_{\text{H}} = \log(\bU) = \frac{1}{2}\log(\bC) = \frac{1}{2} \sum_{i=1}^3 \log(\lambda_i) \bN_i \bN_i^T,
\end{equation}
which is numerically unstable when $\lambda_i \approx 1$. The stable version can be computed by
\begin{equation} \label{EH-stable}
	\bE_{\text{H}} = \frac{1}{2}\log(2\bE + I) = \frac{1}{2} \sum_{i=1}^3 \logp(2\lambda_i^E) \bN_i \bN_i^T.
\end{equation}
We can use a similar approach to define the stable form of the spatial logarithmic in current configuration
\begin{equation} \label{eH-unstable}
	\be_{\text{H}} = \log(\bv) = \frac{1}{2}\log(\bb)
\end{equation}
by using the $\be$'s eigenvalues as
\begin{equation} \label{eH-stable}
	\be_{\text{H}} = \frac{1}{2}\log(2\be + I) = \frac{1}{2} \sum_{i=1}^3 \logp(2\lambda_i^e) \bn_i \bn_i^T.
\end{equation}
With strain computed in a stable way, constitutive models based on the Hencky strain can stably evaluated using the principles in the prior sections.

\section{Algorithmic Differentiation to Derive Material Models}\label{ad}
Deriving some of the material models in \autoref{constitutive} is not a trivial task, and could be tedious and error-prone.
One way to derive complicated material models without needing to manipulate and simplify the intermediate expressions is to use algorithmic (aka. automatic) differentiation (AD).
Starting from the strain energy function in terms of the strain tensor, one could use an AD tool to compute the corresponding stress tensor. 
However, using AD does not \emph{automatically} guarantee stability. Instabilities in the strain energy function will propagate to the derivatives, leading to an unstable evaluation of stress. Our primary goal here is to show how to compute material models using AD, identify instability-inducing terms in the free energy function, and introducing a stable form of the strain energy function. We chose the coupled Neo-Hookean model \eqref{NH-energy} to show the procedure. One can apply these principles to obtain stable representations for other material models.

We start with re-writing equation \eqref{NH-energy} in terms of $\bE$.
\begin{align} \label{NH-energy-E}
	\psi(\bE) &= \frac{\firstlame}{4} \left( J^2 - 1 -2 \log J \right) - \mu \left( \log J - \operatorname{trace}(\bE) \right),
\end{align}
where $J = \sqrt{2 \bE + I}$. As we saw in \eqref{coupled-nh}, equation \eqref{NH-energy-E} is unstable due to the presence of the $J^2 - 1$ and $\log J$ terms. Using \eqref{Jm1}-\eqref{logJ-stable}, we can transform \eqref{NH-energy-E} to
\begin{equation} \label{NH-energy-E1}
	\begin{split}
		\psi(\bE) &= \frac{\firstlame}{4} \left( \mathtt{J_{-1}} \left( \mathtt{J_{-1}} + 2 \right) - 2 \logp \left(\mathtt{J_{-1}} \right) \right) - \mu \left( \logp \left(\mathtt{J_{-1}} \right) - \operatorname{trace}(\bE) \right) \\
		&= \frac{\firstlame}{4} \Big( \underbrace{\mathtt{J_{-1}}^2}_{O(\epsilon^{2})} - 2 \underbrace{\left( \logp  \mathtt{J_{-1}} - \mathtt{J_{-1}} \right)}_{O(\epsilon^{2})} \Big) - \mu \underbrace{\left( \logp \mathtt{J_{-1}} - \operatorname{trace} \bE \right)}_{O(\epsilon^{2})} \, .
	\end{split}
\end{equation}

While \eqref{NH-energy-E1} is more stable than \eqref{NH-energy-E}, it is still unstable when $J \approx 1$ and its derivative results in an unstable formulation for stress due to numerical cancellations in the $\logp  \mathtt{J_{-1}} - \mathtt{J_{-1}}$ and $\logp \left( \mathtt{J_{-1}} \right) - \operatorname{trace}(\bE)$ terms. Looking more closely at these terms, if $\bE$ is of scale $\epsilon$, then $\psi(\bE) \in O(\epsilon^{2})$ and thus we need to avoid subtracting terms such as $\mathtt{J_{-1}}$ and $\trace \bE$, which are $O(\epsilon)$. The first underbrace in \eqref{NH-energy-E1} is fine as is, but the second and third require a reformulation.

For the second underbrace, we define a helper function for computing $\log(1 + x) - x$ that avoids subtracting $O(x)$ terms when computing the $O(x^{2})$ result. Knowing \cite{beebe2017mathematical}
\begin{equation}\label{taylor-log1p}
	\begin{split}
		\logp(x) &= \log(1 + x) \\
		&= 2 \operatorname{artanh} \left( \frac{x}{2 + x} \right) \\
		&= 2 \sum_{n=0}^\infty \frac{1}{2n+1} \left( \frac{x}{2+x} \right)^{2n+1} \, , \, \, \, \big\vert \frac{x}{2+x} \big\vert < 1 \, .
	\end{split}
\end{equation}
and moving $x$ to the left hand side, we have
\begin{equation}\label{taylor-log1p_minus_x}
	\begin{split}
		\logpmx(x) &= \log(1 + x) - x \\
		&= -\frac{x^2}{2 + x} + 2 \sum_{n=1}^\infty \frac{1}{2n+1} \left( \frac{x}{2+x} \right)^{2n+1} , \quad \big\vert \frac{x}{2+x} \big\vert < 1 \, .
	\end{split}
\end{equation}
\autoref{stability-logpmx} studies how many terms are necessary to evaluate $\logpmx$ accurately.

With a stable implementation of $\logpmx$ available, we need only a stable formulation for the third underbrace in \eqref{NH-energy-E1}. Considering a $2\times 2$ strain tensor, we can break down $J^2$ into
\begin{equation} \label{J2}
	\begin{split}
		J^2 &= \Det{\bI + 2 \bE} \\
		&= 1 + 2 \underbrace{\left( E_{1,1} + E_{2,2} \right)}_{\operatorname{trace}(\bE)} + 4 \left( E_{1,1} E_{2,2} - E_{1,2} E_{2,1} \right) \, ,
	\end{split}
\end{equation}
and define a new variable:
\begin{equation} \label{helper-J2m1m2trE}
	\begin{split}
		\mathtt j &= J^2 - 1 - 2 \operatorname{trace}(\bE) \\
		&= 4 \left( E_{1,1} E_{2,2} - E_{1,2} E_{2,1} \right) \, .
	\end{split}
\end{equation}

Using \eqref{helper-J2m1m2trE} extended for $3\times 3$ matrices, we can compute $\mathtt{J_{-1}}$ as
\begin{equation}\label{Jm1-from-J2}
	\mathtt{J_{-1}} = \frac{J^2 - 1}{J + 1} = \frac{\mathtt j + 2 \operatorname{trace}(\bE)}{J + 1} \, .
\end{equation}
Finally, with \eqref{taylor-log1p_minus_x}-\eqref{Jm1-from-J2}, we arrive at a stable representation of the strain energy function \eqref{NH-energy-E} as
\begin{equation}\label{NH-energy-E-stable}
	\psi(\bE) = \frac{\firstlame}{4} \left( \mathtt{J_{-1}}^2 - 2  \logpmx(\mathtt{J_{-1}})  \right) - \mu \big( \logpmx\left(\mathtt j + 2 \trace \bE \right) + \mathtt j \big) \, .
\end{equation}

\begin{figure}[t]
	\centering
	\includegraphics[width=.7\linewidth]{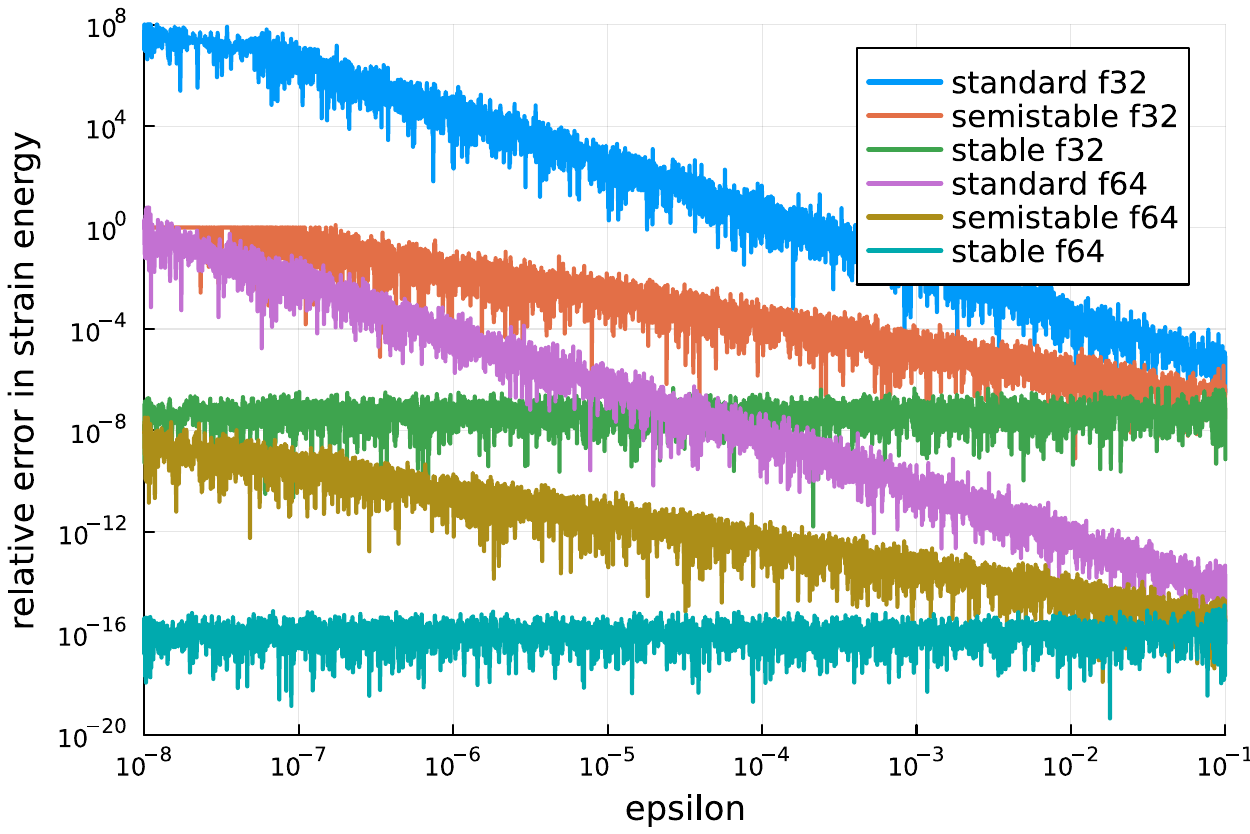}
	\caption{Relative error of computation of Neo-Hookean strain energy in standard \eqref{NH-energy-E}, semistable \eqref{NH-energy-E1} and stable \eqref{NH-energy-E-stable} forms.} \label{fig4:energy}
\end{figure}
\autoref{fig4:energy} shows the relative error of strain energy function for Neo-Hookean model. For the strain of order $10^{-8}$ we lose 16 and 8 digits in the standard \eqref{NH-energy-E} and semistable \eqref{NH-energy-E1} forms, respectively, while the stable form \eqref{NH-energy-E-stable} delivers full accuracy.

\begin{figure}[t]
	\centering
	\includegraphics[width=.7\linewidth]{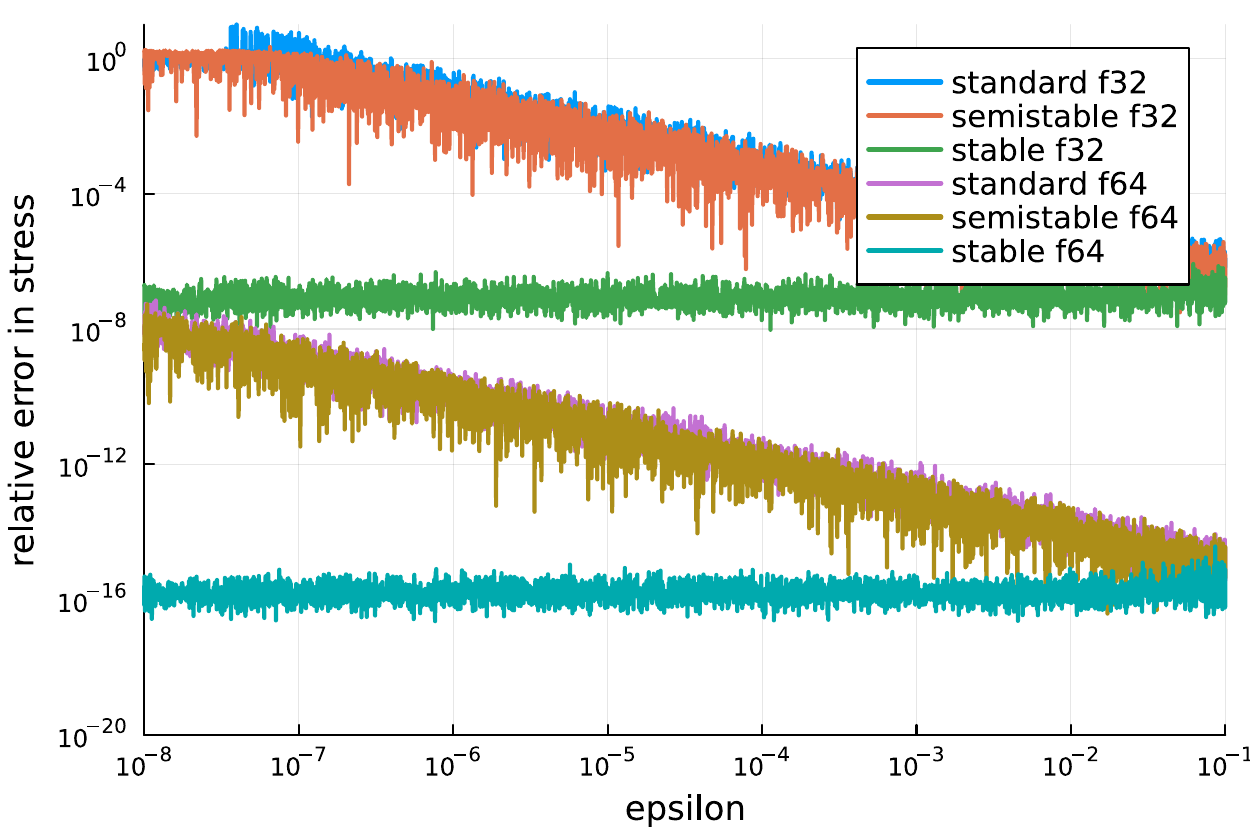}
	\caption{Relative error of computation of stress with AD using the standard, semistable and stable forms of strain energy.} \label{fig5:stress}
\end{figure}

We can now expect AD tools to automatically generate a stable representation of the second Piola-Kirchhoff stress, $\bS = \frac{\partial \psi}{\partial \bE}$, and indeed, \autoref{fig5:stress} shows that direct application of Zygote.jl \cite{Zygote.jl-2018} to \eqref{NH-energy-E-stable} (with $n=6$ in \eqref{taylor-log1p_minus_x}) is stable. Meanwhile, the standard and semistable forms both lose 8 digits for strain at order $10^{-8}$.
\autoref{stress-code} shows sample Julia code to compute $\bS$ via reverse-mode AD.
We note that AD tools that can compute higher derivatives, e.g., by applying forward-mode AD to $\bS(\bE)$, can readily provide ingredients for solvers (such as Newton linearization) and diagnostics, freeing the implementer from tedious coding of higher derivatives or resorting to numerical differentiation.

\begin{lstlisting}[caption={Code for computing stress $\S$ for an arbitrary strain energy $\mathtt{psi}(\bE)$ using the AD tool Zygote.jl.},label=stress-code]
using Zygote
function S(E, psi)
    stress, = Zygote.gradient(psi, E)
    stress
end
\end{lstlisting}

\section{Axial tension/compression test}\label{axial}
We consider an axial tension (or compression) test that illustrates the impact of a stable formulation, and can be used to show when a black-box solver uses an unstable formulation.
Consider a unit elastic block $[0, 1]^{3}$ and apply free slip (symmetry) boundary conditions on the faces $x=0, x=1, y=0, z=0$, with the $x=1$ face translated to $x = 1 + \epsilon$. The reaction forces on the faces $x\in \{0,1\}$ is $\pm E \epsilon + O(\epsilon^{2})$ where $E$ is the Young's modulus.
In our experiment, we will use $E = 2.8, \nu = 0.4$ and $\epsilon = 10^{-12}$ with double precision arithmetic.
The displacement in this setup is a linear function, thus has constant strain and will be exactly represented by any displacement-based conforming finite element discretization at any resolution (including a single cube element).
We consider a Neo-Hookean material with unstable \eqref{NH-S-unstable} and stable \eqref{NH-S-stable} forms, and compute the normal reaction forces on the $x=0$ and $x=1$ faces.
\autoref{tab:axial-newton} demonstrates how rounding error in the unstable form causes Newton to stagnate and leads to inaccurate and unbalanced forces.
This test used a (stable) full-accuracy Jacobian in both tests, so the stagnation is entirely attributable to numerical stability evaluating the residual.
Note that while a formulation cannot be stable without passing this test, it is only one strain tensor, and thus it would be possible for a formulation to be stable for this special case while being unstable generally.

\begin{table}
  \centering
  \begin{tabular}{lrr}
	\toprule
	Iteration & Stable & Unstable \\
	\midrule
	0 & 8.692269873597e-12 & 8.692418103111e-12 \\
	1 & 5.188505917695e-24 & 4.174791583230e-16 \\
	2 & --- & 3.462609857586e-16 \\
	3 & --- & 4.517090534444e-16 \\
	\midrule
	Force $x=1$ & -2.800000000004e-12 & -2.799739663547e-12 \\
	Force $x=0$ & 2.799999999996e-12 & 2.800355568073e-12 \\
	Difference & -8e-24 & 6.15905e-16 \\
    \bottomrule
  \end{tabular}
  \caption{Newton exhibits quadratic convergence for the stable form, but stagnates for the unstable form. The reaction force is accurate to 12 digits and balanced between the axial faces for the stable form, but has only four digits of accuracy and associated force imbalance for the unstable form.}
  \label{tab:axial-newton}
\end{table}

\section{Conclusions}\label{conclusion}
In this paper, we investigated various constitutive formulations for elasticity along with their implementation in finite element software packages. 
Formulations written in terms of the deformation gradient $\bF$ cannot be numerically stable.
Standard formulations have additional instabilities due to the presence of function like $J-1$ and/or $\log J $, which are unstable when $J \approx 1$, as well as terms like $\bI - \bC^{-1}$, which similarly experience cancellation for small strain $\bE \approx 0$.
In general, the standard computation for a strain of order $10^{-m}$ will result in $m$ digits lost in the computed stress tensor and $2m$ lost evaluating the strain energy function.
We proposed equivalent stable formulations, all of which achieve relative accuracy $O(\emachine)$.
These new formulations make use of the displacement gradient $\bH$ to define a strain tensor without loss of significance, compute $J-1$ in a stable way, and avoid cancelation computing shear stress terms.
In addition to coupled and decoupled forms of Neo-Hookean, Mooney-Rivlin, and Ogden models in initial and current configuration, we showed that one can achieve a stable formulation using algorithmic differentiation (AD) if the forward model (strain energy) is stable.
We also showed a stable evaluation of Hencky strain, which is important for developing stable representations for inelastic models at finite strain.

With single precision in the standard formulation, the first digit of stress is incorrect for strains of order $10^{-7}$, while the new stable formulations get all 7 digits correct at all strain levels.
The stable formulations open the door for hyperelastic simulation using single or mixed precision \cite{abdelfattah2021survey}, thereby improving performance and reducing hardware and energy cost without compromising accuracy.
Moreover, stable formulations are necessary to run efficiently on hardware that does not support double precision, such as GPUs, tensor cores, and embedded devices.

\section*{Data availability}
The software/data that support the findings of this study are openly available in Zenodo at \url{https://doi.org/10.5281/zenodo.11113499} \cite{shakeri2024stabledemo}.

\bibliography{refs}%

\appendix
\section{Numerical Stability Evaluation}\label{num-eval-stability}

In order to compare stability of different implementations of functions of the displacement gradient $\bH = \frac{\partial \bm u}{\partial \bm X}$, we start with by sampling $\bH = \hat{\bH}/\lVert \hat{\bH} \rVert_{F}$, in which $\hat{\bH}$ is the absolute value of a standard normal distribution (\texttt{Hhat = randn(3, 3)} in Julia) and then plot relative error of each function $f(\epsilon \bH)$ as in \autoref{relative-error-code}, where $\epsilon \in (10^{-8}, 10^{-1})$ to cover a range from small to large strain. Julia's \texttt{big} converts the input to arbitrary precision (default gives $\emachine < 10^{-77}$) and further operations retain that arbitrary precision. For the range of $\epsilon$ considered, the \texttt{big} arithmetic can be considered an exact reference and then calculate the relative error for single and double precision i.e., \texttt{repr = Float32} and \texttt{repr = Float64}.

\begin{lstlisting}[caption={Julia code for computing relative error of function $f(\epsilon\bH)$},label=relative-error-code]
function rel_error(eps, f, repr)
    H_hat = randn(3,3)
    H = H_hat / norm(H_hat)
    
    ref = f(big.(eps*H)) # arbitrary precision
    norm(f(repr.(eps*H)) - ref) / norm(ref)
end
\end{lstlisting}

We start with the $J-1$ term which appears in all hyperelastic models and compare it with its stable form $\mathtt{J_{-1}}$ \eqref{Jm1} as defined in Julia in \autoref{Jm1-code} and their relative errors are shown in figure \eqref{fig1:Jm1}. As expected, the stable computation $\mathtt{J_{-1}}$ has a relative accuracy of order $\emachine$ for single and double precision and $J-1$ loses accuracy as $\bH$ decreases. In fact, for $u_{i,j}$ of order $10^{-8}$ we can trust no digits in single precision and we lost half of the digits in double precision, respectively.

\filbreak
\begin{lstlisting}[caption={Code for computing $J-1$ and $\mathtt{J_{-1}}$ },label=Jm1-code]
function Jm1_unstable(H)
    F = I + H
    J = det(F)
    J - 1
end

function Jm1_stable(H)
    det_H = det(H)
    A1 = H[1,1]*H[2,2] + H[1,1]*H[3,3] + H[2,2]*H[3,3]
    A2 = H[1,2]*H[2,1] + H[1,3]*H[3,1] + H[2,3]*H[3,2]
    # Compute J-1
    det_H + tr(H) + A1 - A2
end
\end{lstlisting}

To assess stability of constitutive models, we implement the standard and unstable expressions for the appropriate stress $\bS(\bH)$ or $\btau(\bH)$, internally making use of Green-Lagrange or Green-Euler strains computed by stable means \eqref{eq:green-lagrange-stable}, and measure the relative error via \autoref{relative-error-code}, yielding figures like \autoref{fig2:Ogden} and \autoref{fig5:stress}.

\section{The \texttt{log1p} family of functions}\label{appendix:log1p}
The $\logp$ function is available, for example, in C99 and 4.3BSD \texttt{math.h}, Python \texttt{math.log1p}, Julia \texttt{log1p}, Rust \texttt{f64::ln\_1p}, with similar interfaces for $\expm$.
Although these functions are specified by IEEE 754-2008 (named \texttt{logp1} and $\expm$), they are conspicuously missing in Fortran as of the Fortran 2023 standard.
Structural mechanics software written in Fortran may access the routines from C's \texttt{libm} using an interface binding as in \autoref{log1p-code}, by using identities in terms of the hyperbolic tangent and its inverse,
\begin{align*}
	\logp(x) &= \log(1+x) = 2 \operatorname{artanh}\left(\frac{x}{2 + x} \right) & \expm(x) &= \exp(x) - 1 = \frac{2\tanh(x/2)}{1 - \tanh(x/2)},
\end{align*}
or approximated using a series expansion such as \eqref{taylor-log1p}.
We have observed that compilers/math libraries often do not provide vectorized support for functions such as $\logp$ and \texttt{atanh}, in which case a series expansion valid to sufficient accuracy across the range of physically-realizable inputs offers a performance benefit.

\begin{lstlisting}[caption={Sample Fortran-2003 interface binding a Fortran-callable function named \texttt{logp1} (the IEEE 754-2008 name, and likely name in a future Fortran standard) to \texttt{log1p} from the C math library \texttt{libm}.},label=log1p-code]
interface
   function logp1(x) bind(c, name="log1p")
   use iso_c_binding, only: c_double
   real(c_double), intent(in), value :: x
   real(c_double) :: logp1
   end function
end interface
\end{lstlisting}

\section{Accurate evaluation of $\logpmx$}\label{stability-logpmx}
\begin{figure}[t]
	\centering
	\includegraphics[width=.7\linewidth]{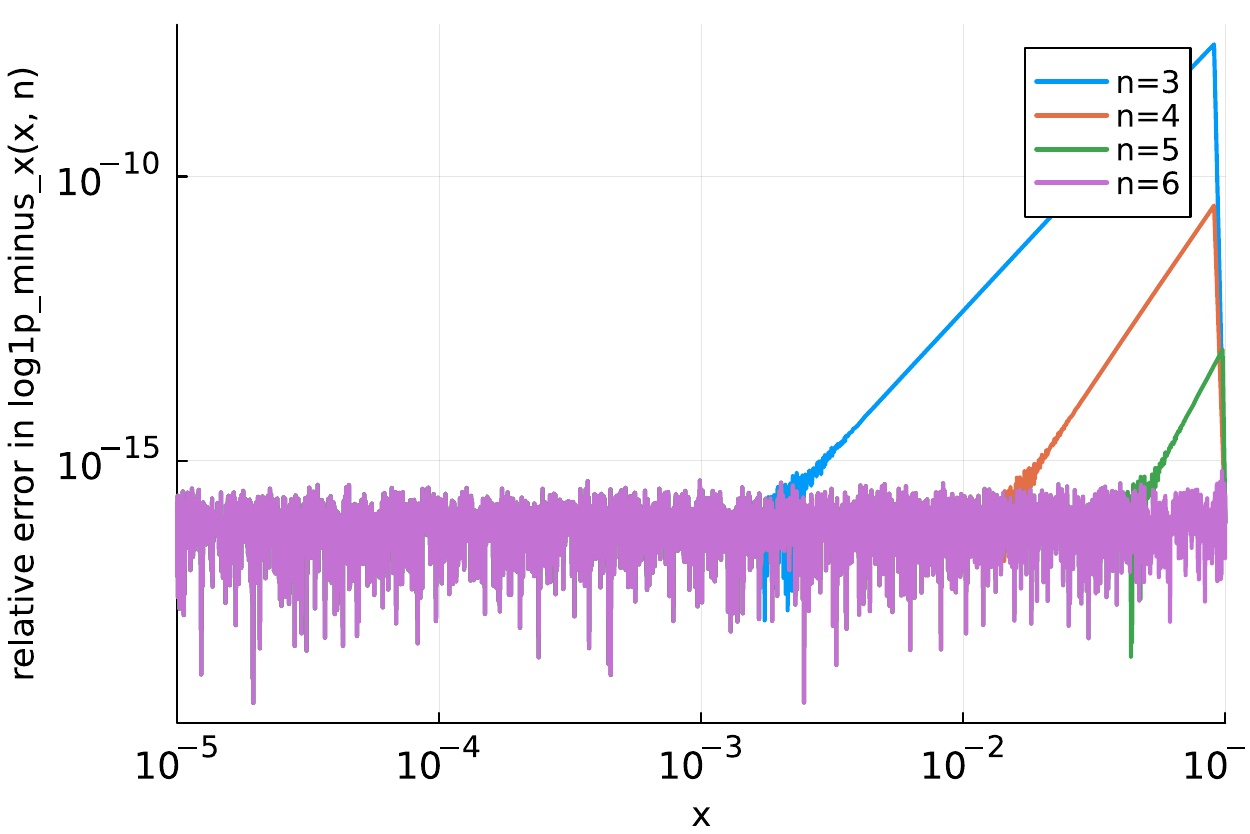}
	\caption{Relative error of computation of $\logpmx()$ for different expansions of the Taylor series.} \label{fig3:taylor-n}
\end{figure}

We require efficient and accurate evaluation of the function $\logpmx(x) = \log(1 + x) - x$ to ensure accuracy of the AD-computed stresses in \autoref{ad}.
\autoref{fig3:taylor-n} demonstrates that $n=6$ in \eqref{taylor-log1p_minus_x} is sufficient to provide $O(\emachine)$ accuracy evaluating $\logpmx(x)$ for $\epsilon \in (10^{-5}, 10^{-1})$.
We believe it would be fruitful to develop a uniformly accurate algorithm for $\logpmx$ and include it in numerical libraries.

\end{document}